\documentclass[11pt]{interact}
\usepackage{epstopdf,comment,ulem}
\usepackage{subfigure}
\usepackage{color,soul}

\usepackage{amsmath,amssymb,amsfonts}
\usepackage{pifont}

\makeatletter
\newcommand{\centeredcaption}[2][]{%
  \begingroup
  \renewcommand{\@makecaption}[2]{\centering\small\textbf{##1: }##2\par}
  \caption[#1]{#2}%
  \endgroup
}
\makeatother

\usepackage{graphicx}
\usepackage{float}
\usepackage[dvipsnames]{xcolor}

\usepackage[numbers,sort]{natbib}

\bibpunct[, ]{(}{)}{;}{n}{}{,}

\theoremstyle{plain}
\newtheorem{theorem}{Theorem}[section]

\newtheorem{proposition}[theorem]{Proposition}

\theoremstyle{definition}

\theoremstyle{remark}
\newtheorem{remark}{Remark}

\begin{document}
\setcitestyle{square}

\articletype{RESEARCH ARTICLE}

\title{Optimization of a lattice spring model with elastoplastic conducting springs: A case study}

\author{
\name{Sakshi Malhotra\textsuperscript{a},Yang Jiao\textsuperscript{b}\thanks{*YJ is supported by  grant NSF CMMI-1916878} and Oleg Makarenkov\textsuperscript{c}\thanks{ OM is supported by grant NSF CMMI-1916876.}\\
\affil{\textsuperscript{a}Sakshi Malhotra is with Department of Mathematical Sciences,
        University of Texas at Dallas, Richardson, TX 75080, USA {\tt\small Sakshi.Malhotra@utdallas.edu;}} \textsuperscript{b}Yang Jiao is with the Material Sciences and Engineering Program, School for Engineering of Matter, Transport and Energy, Arizona State University, Tempe, AZ 85281, USA {\tt\small yang.jiao.2@asu.edu;}}\textsuperscript{c}Oleg Makarenkov is with the Department of Mathematical Sciences, University of Texas at Dallas, Richardson, TX 75080, USA
        {\tt\small makarenkov@utdallas.edu}}

\maketitle

\begin{abstract}
We consider a simple lattice spring model in which every spring is elastoplastic and is capable to conduct current.
The elasticity bounds of spring $i$ are taken as $[-c_i,c_i]$ and the resistance of spring $i$ is taken as $1/c_i$,
which allows us to compute the resistance of the system.
The model is further subjected to a gradual stretching and, due to plasticity, the response force increases until 
a certain terminal value. 
We demonstrate that the recently developed sweeping process theory can be used to optimize the interplay between the terminal response force
and the resistance on a physical domain of parameters $c_i.$ The proposed methodology can be used by practitioners for the design of multi-functional materials as an alternative to topological optimization. 
\end{abstract}

\begin{keywords}
component, formatting, style, styling, insert
\end{keywords}

\section{Introduction}

Multi-functionality, e.g. when the stiffness and electrical conductivity of a composite are simultaneously optimized, is a key principle for novel 
lightweight 3D printed network materials (also called metamaterials). Such materials are highly desirable for a wide spectrum of applications from aerospace engineering to soft robotics and flexible electronics \cite{Wang,mixed2,Yang-mixed}. Despite of the fact that, in the presence of mechanical loading, the reliable estimation of plastic behavior of the material requires accounting for plasticity of individual bonds of material's microstructure \cite{blechman,survey-printing}, current systematic theories for optimization of multi-functional properties focus on just elastic materials. The goal of the paper is to examine the opportunities that the recently developed theory of sweeping processes \cite{ESAIM,SICON} opens for optimization of multi-functional materials with plasticity.

\vskip0.2cm

\noindent To model multi-functional materials with plasticity, this paper proposes to employ the so-called lattice spring model where the sweeping process theory \cite{Moreau,ESAIM} recently provided algebraic formulas \cite{SICON} for the response force and where the springs can simultaneously be viewed as conductors. Efficient methods that map actual materials to lattice spring models are developed in e.g. Chen et al \cite{Chen} and Kale and Ostoja-Starzewski \cite{Kale}. 

\vskip0.2cm



\begin{figure}[h]
\centering
\includegraphics[scale=0.6]{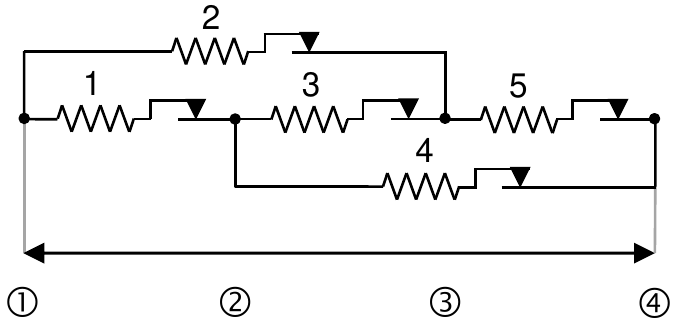}
\caption{A system of 5 springs on 4 nodes with a displacement-controlled loading that locks the distance between nodes \ding{172} and \ding{176}.} \label{fig1}
\end{figure}

\noindent The analysis in this paper concerns the lattice spring model of Fig.~\ref{fig1}, because this is the minimal irreducible model that requires a 3D (i.e. non-planar) sweeping process framework (Rachinskiy \cite{rach}). Therefore, we expect that all interesting non-trivial optimization phenomena are to be captured by the model of Fig.~\ref{fig1}. In other words, we expect that model of Fig.~\ref{fig1} is an efficient benchmark to assess the capabilities of the sweeping process framework in the optimization of most typical multi-functional properties of  lattice-spring models of general form. Indeed, using the model of Fig.~\ref{fig1} we conjectured a number of qualitative properties of universe nature (see Conclusions section). On the other hand, we included quotes from Wolfram Mathematica code and a computational guide (Appendix~\ref{appB}) that allow to generalize simulations of this paper to  arbitrary lattice-spring models with 1-dimensional nodes (as long as the resistance of the network is computable via the known methods of e.g. \cite{c,resistance}). Further extension of the sweeping process framework to the case of nodes of arbitrary dimension and hyperuniform network materials is available in \cite{Gud} whose formulas can also be used for optimization of multi-functional properties along the strategy of the current paper.

\vskip0.2cm

\noindent In what follows, $c_i$ denotes the elastic limit of spring $c_i$, so that the relaxed length of the spring changes when the stress attempts to exceed the interval $[-c_i,c_i].$ For lattice spring models coming from fiber-reinforced composites, the sum $c_1+c_2+c_3+c_4+c_5$ 
can be viewed as the fabrication cost of the material
because the cost of such composites is mainly determined by fabrication of the nano-carbon fibers \cite{mater01}. When the displacement-controlled loading (applied to nodes \ding{172} and \ding{175}) stretches the system gradually, the response force of the system will increase up to a certain terminal value $F$ that quantifies the strength of the material. Viewing Fig.~\ref{fig1} as an electric circuit, and using $R$ and $G$ to denote the resistance and conductance between nodes \ding{172} and \ding{175}, the present paper addresses optimization of the functionals 
\begin{equation}\label{functionals1}
  F_R=k_1 F+k_2 R,\quad F_G=k_1 F+k_2 G,\quad {\rm and}\quad   C=c_1+c_2+c_3+c_4+c_5.
\end{equation}


\noindent Using Wolfram Mathematica we make a variety of qualitative observations (listed in Conclusions section) that could motivate a good deal of new theoretical research in optimization theory. Some of these observations are uniqueness of optimal values of parameters $(c_1,...,c_5)$ in the optimization of $F_R$ and $F_G$, non-uniqueness of $(c_1,...,c_5)$ in the optimization of $C$,  existence of linear thresholds in the space of parameters that separate topologically different patterns of optimal $(c_1,...,c_5).$

\vskip0.2cm

\noindent The paper is organized as follows. In the next section of the paper we introduce the ingredients $F$, $R$, and $G$ used in the definition of the functionals $F_R$ and $F_G$. In particular, we introduce two domains (\ref{FF2}) and (\ref{FF4}) for the parameters $c_1,...,c_5$ that this paper restricts the analysis to (taken from Gudoshnikov et al \citep{SICON}). Section~\ref{methods} discusses the numerical methods used for optimization of the functionals $F_R$, $F_G$, $C$, {\color{black} along with computational complexity of the method of our choice.} Section~\ref{costsection} is devoted to maximization of the functionals $F_R$ and $F_G$ with a constraint on the cost $C.$ In this section we first determine the maximal value $F_{max}$ of $F$ (Section~\ref{maxF}) that can be reached by the system of Fig.~\ref{fig1} with the chosen constraint on the cost and then maximize $F_R$ and $F_G$ (Sections~\ref{SectionMR} and \ref{Sec33} respectively) restricting on the domain of $c_1,...,c_5$ for which the value of $F$ doesn't go below $F_{max}$ by more than $25\%$ (constraint (\ref{optF})). In other words, we forbid improved multi-functional properties of a lattice-spring system to compromise the strength of the system by more than $25\%.$ Minimization of the cost $C$ under constraint on $F_R$ and $F_G$ (from below) is addressed in Sections~\ref{FR05} and \ref{FG05}. Sections~\ref{FR05} and \ref{FG05} are restricted on the domain of $c_1,...,c_5$ for which the value of $F$ doesn't go below $F_{max}$ too. The conclusions from numerical simulations of Sections~\ref{costsection} and \ref{Sec4} are formulated as Propositions~\ref{prop32}, \ref{prop33}, \ref{prop41}, \ref{prop42}.
Additionally, Section~\ref{maxF} quotes the Mathematica's command we use, which clarifies that all numerical computations can be executed analytically by hand over the methods of linear programming. 
 The conclusion section interprets the numerical results of Sections~\ref{costsection} and \ref{FR05} from more general theoretic perspectives and draws a list of analytic questions inspired by Propositions~\ref{prop32}, \ref{prop33}, \ref{prop41}, \ref{prop42}.
  Two appendixes that compute the resistance of the circuit of Fig.~\ref{fig1} (Appendix~\ref{appA}) and outline the sweeping process framework of \cite{SICON} based conclusions about plastic regimes in elastoplastic system of Fig.~\ref{fig1} (Appendix~\ref{appB}) conclude the paper. Wolfram Mathematica Notebook is attached as supplementary material.

\section{Algebraic formulas for the optimization functionals}

\noindent {\bf Resistance of the System of Springs.} Assuming that the elastic limits of the 5 springs of Fig.~\ref{fig1} are given by $[-c_i,c_i]$, $i\in\overline{1,5}$, the resistance $R_i$ of spring $i$ will be taken as
\begin{equation}
   R_i=1/c_i.
\label{eq_Ri}  
\end{equation}
The relation (\ref{eq_Ri}) is  inspired by typical design mechanisms used in composite materials. For example, one possible realization of the spring can be achieved by fiber-reinforced polymer matrix composites \cite{mater01}. In this system, both the mechanical strength (i.e. $c_i$) and conductivity (i.e. $1/R_i$) are determined by the fiber phase and thus, are positively correlated \cite{mater02}.
  With the values of resistances of springs given by $R_{i}$, the resultant resistance between nodes \ding{172} and \ding{175} computes as (see Appendix~\ref{appA})
\begin{equation}
    R=\dfrac{c_3c_4+c_2(c_3+c_4)+c_3c_5+c_4c_5+c_1(c_2+c_3+c_5)}{c_1c_4c_5+c_2c_4c_5+c_1c_2(c_4+c_5)+c_2c_3(c_4+c_5)}
    \label{eq_R}
 \end{equation}

\noindent and conductance is computed as
\begin{equation}
  G=1/R.
   \label{eq_Cd}
\end{equation}

\noindent {\bf Terminal Response Force.}  According to the sweeping process framework by Gudoshnikov et al \cite{SICON} (see Appendix~\ref{appB} for details), gradual stretching of the lattice of Fig.~\ref{fig1} via the displacement-controlled loading that is applied between nodes \ding{172} and \ding{175} increases the response force of the lattice up to a certain terminal value whose formula takes different forms depending on the values of $c_1,...,c_5$. These formulas can be summarized as follows
\begin{eqnarray}
  F&=&c_1+c_3+c_5, \label{FF1} \\
  &&{\rm if}\ \ -c_2\le c_3+c_5\le c_2, \ \ -c_4\le c_3+c_1\le c_4, \label{FF2} \\
  F&=& c_2+c_3+c_4, \label{FF3} \\
   &&{\rm if}\ \ -c_1\le c_3+c_4\le c_1,  \ \ -c_5\le c_2+c_3\le c_5, \label{FF4}
\end{eqnarray}
where we restrict ourselves to parameters that ensure the greatest number of springs reaches the plastic mode. The later property corresponds to the most uniform distribution of plastic deformation across the lattice that reduces the risk of crack initialization. 
Specifically, if inequalities (\ref{FF2}) hold in a strict sense then springs 1, 3, 5 reach plastic deformation as the magnitude of the displacement-controlled loading crosses a threshold (i.e. {\it springs 1, 3, 5 get to the plastic mode terminally}, 
see \cite{SICON} for explicit formulas), while if condition (\ref{FF4}) holds in a strict sense then those are springs 2, 3, 4 that get to plastic mode terminally\footnote{Note, the non-strict inequalities in (8.4) and (8.5) of \cite{SICON} is a typo. There should be strict inequalities because (8.4) and (8.5) come from (7.10) whose inequalities are strict.}. The inequalities  (\ref{FF2}) and (\ref{FF4}) will be referred to as {\it feasibility conditions.} In what follows, we use feasibility conditions as optimization constraints, so that we will usually get optimal $c_i$ on the boundary of the domains (\ref{FF2}) and (\ref{FF4}) and  strict inequalities in (\ref{FF2}) and (\ref{FF4}) won't hold for our optimal $c_i$. Arbitrary small perturbation will be needed to move such optimal $c_i$ to the interior of the domain (\ref{FF2}) or (\ref{FF4}) to make the plastic behavior of the model predictable according to \cite{SICON}.

\vskip0.2cm

\begin{remark}\label{uniquenessremark} We note that if $(c_1,c_2,c_3,c_4,c_5)$ satisfies (\ref{FF2}) then 
\begin{equation}\label{crelation}
(\tilde c_1,\tilde c_2,\tilde c_3,\tilde c_4,\tilde c_5)=(c_2,c_1,c_3,c_5,c_4)
\end{equation} satisfies (\ref{FF4}) and $R,F$ computed for $(c_1,c_2,c_3,c_4,c_5)$ coincides with $(\tilde R,\tilde F)$ computed for 
$(\tilde c_1,\tilde c_2,\tilde c_3,\tilde c_4,\tilde c_5)$. Therefore, solving an optimization problem involving $F$ and $R$ on the domain (\ref{FF2}) is equivalent to solving same optimization problem on the domain (\ref{FF4}). We will, however, run Wolfram Mathematica code for both problems to assess uniquness/nonuniqueness of  optimal values $(c_1,...,c_5)$. In other words, if Mathematica code executed for domains (\ref{FF2}) and (\ref{FF4}) returns us $(c_1,...,c_5)$ and $(\tilde c_1,...,\tilde c_5)$ that are related over (\ref{crelation}), then we conclude that the optimal value $(c_1,...,c_5)$ is  unique in the domain (\ref{FF2}) and that the optimal value $(\tilde c_1,...,\tilde c_5)$ is  unique in the domain (\ref{FF4}).
\end{remark}





\section{The methods used for numerical simulations}\label{methods}

\noindent We analyze the results for the optimization problem based on two different optimization methods for constrained global optimization available within Wolfram Mathematica. We use NMaximize and NMinimize which are global optimization algorithms for solving problems numerically. Within those we have other optimization methods, the ones that we consider for our problem are "RandomSearch" and "Differential Evolution". Differential evolution is a simple stochastic function minimizer. It starts \cite{ref2, diffevo} with a population of candidate solutions randomly distributed in the search space. At each iteration, new candidate solutions are generated by combining and mutating existing solutions according to certain rules. These new solutions are then evaluated, and the best ones are retained for the next generation. Random search on the other hand works \cite{ref2} by generating a population of random starting points and uses a local optimization method from each of the starting points to converge to a local minimum. {\color{black} When the results from the two methods do not coincide, Differential Evolution was the one showing the correct extremum except of a negligible number of exceptions that we furter verified using Phyton SLSQP algorithm. We therefore conclude that Differential Evolution is a reliable method for the optimization problem considered in this paper.} 

\vskip0.2cm

\noindent {According to \cite{DEcomplexity}, the computational complexity of the method of Differential Evolution for large networks is determined by the dimensions of the search space, the population size, and the number of iterations. Specifically, the general computational complexity for Differential Evolution is given by \cite{DEcomplexity} 
\begin{equation}\label{complexity}
     \mathcal{O} (n\cdot N_{p}\cdot G_{max})
\end{equation}
under assumption that the computational complexity of evaluating the objective function grows at most linearly with the dimension of the search space (which is the case for the objective functions of the present paper because $C$ is linear by construction, $F$ is linear by \cite[formula (3.23)]{ESAIM}, and computational time for $R$ scales close to linear if using  the algorithm \cite{resistancecomplexity}, see \cite[table~I]{resistancecomplexity}). Formula (\ref{complexity}) says that the number of operations required by Differential Evolution is proportional to the product of the problem's dimensionality ($n$), the population size ($N_p$), and the number of iterations ($G_{max}$). The population size $N_{p}$ is usually a fixed constant between $5n$ and $10n$ \cite{DEcomplexity2}. The number of iterations $G_{max}$ depends on the desired accuracy and can be viewed as another fixed constant independent of $n$ and $N_p$. Therefore, if the population size is proportional to $n$, then computational complexity of the method of Differential Evolution scales quadratically with respect to the dimension of the problem. 
}

\section{Maximization of multi-functional properties under constraint on the cost of the material}\label{costsection}

\noindent In this section we search for the parameters $c_i$ that optimize $M_R$ and $M_G$ under assumption that 
\begin{equation}\label{cost1}
  c_1+c_2+c_3+c_4+c_5 \leq 2.
\end{equation}

\subsection{Maximization of strength in the absence of electrical properties}\label{maxF}

We first determine the maximal possible strength of the material (i.e. the maximal possible $F$) in the absence of conductance and resistance because we want that addition of conductivity and resistance does not compromise the strength by more than 25\%, which is desirable in, e.g., soft robotic applications \cite{mater01, Yang-mixed, Xu_JAP,strength-yang}. 


\vskip0.2cm


\noindent From definition (\ref{FF1})-(\ref{FF4}) we conclude that maximization of the response force splits into two cases:

\vskip0.2cm

\noindent{\it Case 1:} (\ref{FF2}) holds. Using Wolfram Mathematica with the command 

\vskip0.1cm
\begin{itemize}\item[]  Maximize $[ \{ c_1+c_3+c_5 , c_1\geq 0, c_2\geq 0, c_3\geq 0, c_4\geq 0, c_5\geq 0, -c_3-c_2-c_5 \leq 0,$ \\
  $c_3+c_5-c_2 \leq 0, -c_1-c_3-c_4 \leq 0 , c_1+c_3-c_4 \leq 0 ,c_1+c_2+c_3+c_4+c_5 \leq 2 \} , \{ c_1,c_2,c_3,c_4,c_5 \} ] $ 
\end{itemize}

\vskip0.1cm

\noindent to maximize $F$ given by (\ref{FF1}) under constraints (\ref{FF2}) and (\ref{cost1}) leads to
\begin{equation}\label{Fmax}
  F=1,
\end{equation}
which value is attained at $c_2=c_3=c_5=0$, $c_1=c_4=1$

\vskip0.2cm

\noindent{\it Case 2:} (\ref{FF4}) holds. Using Wolfram Mathematica (with a command analogous to Case~1) to maximize $F$ given by (\ref{FF3}) under constraints (\ref{FF4}) and (\ref{cost1}) leads to the same value (\ref{Fmax}) of $F$, which is now attained at $c_1=c_3=c_4=0$, $c_2=c_5=1$ (mirroring Case 1 as expected).

\vskip0.2cm

\noindent The results of Case 1 and Case 2 can be concluded as follows.

\begin{figure}[h]
\centering
\includegraphics[scale=0.6]{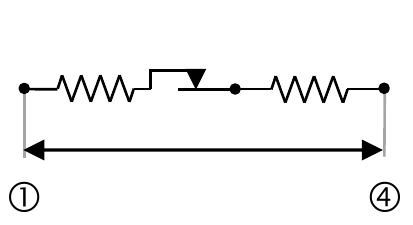}
\caption{The model of Fig.~\ref{fig1} with springs 2,3,5 removed (i.e. with $c_2=c_3=c_5=0$). The arrow in the left spring indicates that the left spring will deform plastically as the displacement-controlled loading (double-sided arrow) increases.} \label{figF}
\end{figure}

\begin{proposition}  In order to obtain the maximal strength of the network of Fig.~\ref{fig1} keeping the cost of the network under a fixed value, only 2 springs connected in serial as shown in Fig.~\ref{figF} are needed. 
\end{proposition}


\noindent Based on the value (\ref{Fmax}) we impose the following additional constraint 
     \begin{equation}\label{optF}
         F \geq 0.75,
     \end{equation}
saying that we don't want to   reduce the strength by more than 25\% when optimizing conductance and resistance.


\subsection{Maximization of the strength-resistance functional $F_R$}\label{SectionMR}

 In this section we use Wolfram Mathematica to find the maximal value of the functional $F_R$ (given by (\ref{functionals1})) on set of parameters $c_1,...,c_5$ satisfying  constraints (\ref{cost1}) and (\ref{optF}) along with either (\ref{FF2}) or (\ref{FF4}). This optimization problem is solved for each fixed pair of coefficients $k_1,k_2$ varying from 0.1 to 1 with the step 0.1. We use two commands of the code: 
\begin{itemize}
\raggedright

\item []{\color{black} Table $[ \{ K_1, K_2, c_1, c_2, c_3, c_4, c_5, F_7, c_1+c_2+c_3+c_4+c_5, \#  \} \, / \, . \, \# \, 2 \, \& \, @ \,@ \,$ \\
NMinimize $[ \{  K_1*F_7 + K_2* R[c_1,c_2,c_3,c_4,c_5], c_1 \geq 0 , c_2\geq 0, c_3 \geq 0,c_4 \geq 0,c_5 \geq 0 , c_1+c_2+c_3+c_4+c_5 \leq 2 , F_7 \geq 0.75 , -c_3-c_2-c_5 \leq 0 , c_3+c_5-c_2 \leq 0 , -c_3 -c_1-c_4 \leq 0 , c
_3+c_1-c_4 \leq 0 \}, \{ c_1 , c_2 ,c_3 ,c_4 , c_5\} , \text{Method} \rightarrow \{ "\text{RandomSearch}"\}, "\text{AccuracyGoal}" \rightarrow 6  ], \{ K_1 , 0.1 , 1, 0.1 \} , \{ K_2 , 0.1 , 1, 0.1 \} 
 ] \thicksim \text{Flatten} \thicksim 1; $} 
 \item []{\color{black} Table $[ \{ K_1, K_2, c_1, c_2, c_3, c_4, c_5, F_7, c_1+c_2+c_3+c_4+c_5, \#  \} \, / \, . \, \# \, 2 \, \& \, @ \,@ \,$ \\
NMinimize $[ \{  K_1*F_7 + K_2* R[c_1,c_2,c_3,c_4,c_5], c_1 \geq 0 , c_2\geq 0, c_3 \geq 0,c_4 \geq 0,c_5 \geq 0 , c_1+c_2+c_3+c_4+c_5 \leq 2 , F_7 \geq 0.75 , -c_3-c_2-c_5 \leq 0 , c_3+c_5-c_2 \leq 0 , -c_3 -c_1-c_4 \leq 0 , c
_3+c_1-c_4 \leq 0 \}, \{ c_1 , c_2 ,c_3 ,c_4 , c_5\} , \text{Method} \rightarrow \{ "\text{DifferentialEvolution}"\}, "\text{AccuracyGoal}" \rightarrow 6  ], \{ K_1 , 0.1 , 1, 0.1 \} , \{ K_2 , 0.1 , 1, 0.1 \} 
 ] \thicksim \text{Flatten} \thicksim 1; $}
\end{itemize}
and choose the best result for each of the values of $k_1$ and $k_2$.

\vskip0.2cm

\begin{figure}[H]
\centering
\includegraphics[scale=0.65]{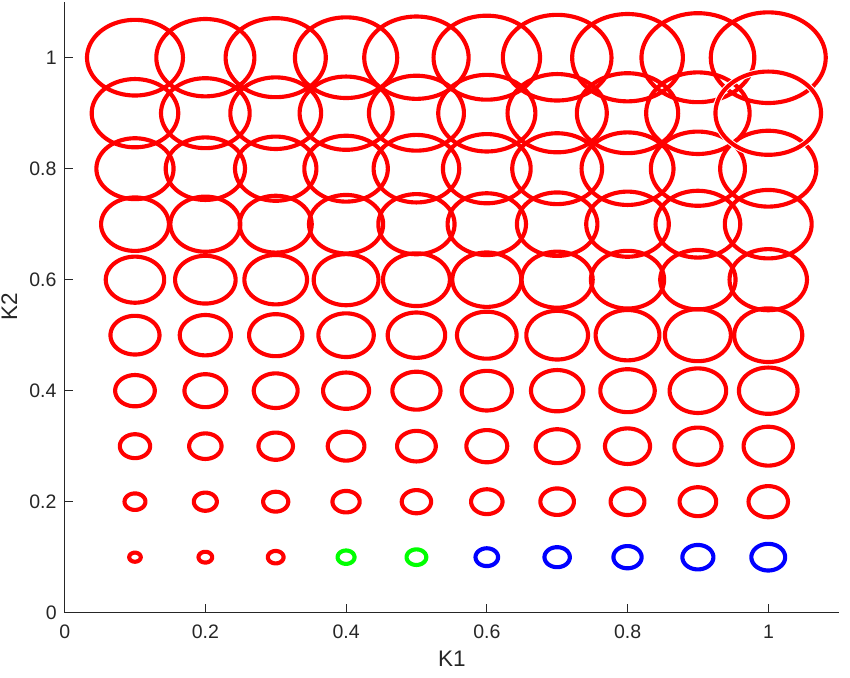}
\caption{Maximization of the functional $F_R$. The horizontal and vertical axis stay for $k_1$ and $k_2$ respectively. The radius of the disc centered at $(k_1,k_2)$ refers to the maximum of $F_R$ for the parameters  $(k_1,k_2)$. The circles in red denote the points where the response force $F$ takes the value 0.75 and resistance of the system $R$ takes the value 3.33332. The circles in blue denote the points where $F=1$ and $R=2$. The circles in green are the cases where all the optimal values of $c_i^s$ are non zero.  } \label{fig3a}
\end{figure}

\begin{figure}[H]
\centering
\includegraphics[scale=0.75]{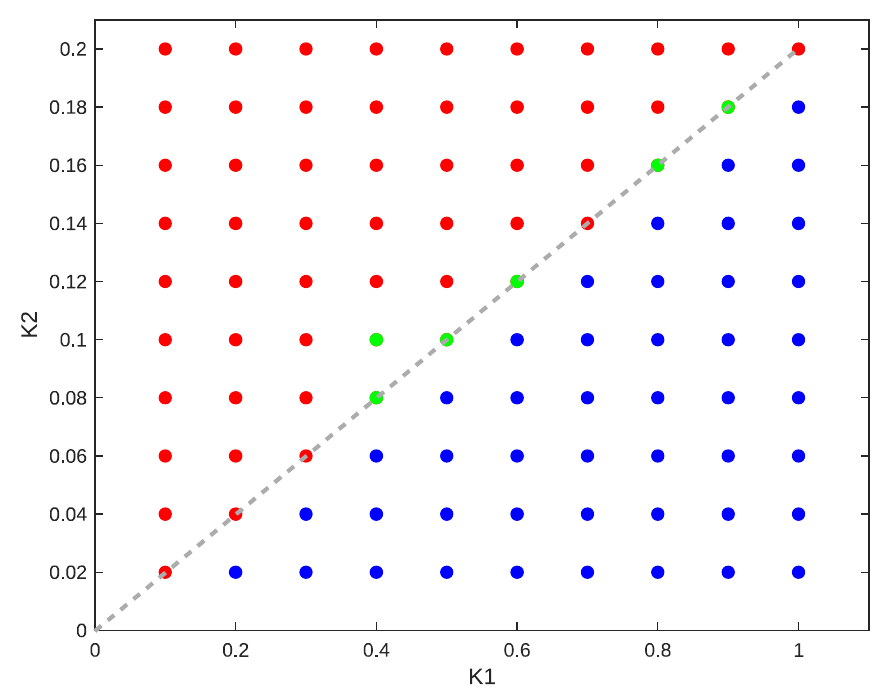}
\caption{Fig.~\ref{fig3a} zoomed in.} \label{fig3a-zoom}
\end{figure}

\noindent We conclude that for each fixed $k_1$ and $k_2$, the maximal value of $F_R$ and the corresponding values of $F$ and $R$ do not depend on whether constraint (\ref{FF2}) or (\ref{FF4}) is used. These values are presented in Fig.~\ref{fig3a}. Increase of the maximal value of $F_R$ when $k_1$ and $k_2$ increase is expected because of the structure of the functional $F_R$. However, the switch from $(F,R)=(0.75,3.33)$ to $(F,R)=(1,2)$  when $(k_1,k_2)$ crosses a threshold (solid gray line in Fig.~\ref{fig3a-zoom}) is an interesting discovery.

\vskip0.2cm

\noindent The optimal values of $c_1,...,c_5$ at which the maximum of $F_R$ is attained turned out to be strictly related to whether $(F,R)=(0.75,3.33)$ or $(F,R)=(1,2)$   holds. 

\vskip0.2cm

\noindent For the points $(k_1,k_2)$ where $(F,R)=(0.75,3.33)$ {\color{black}(i.e. above the threshold of Fig.~\ref{fig3a-zoom})}, we have $(c_1,c_2,c_3,c_4,c_5)=(0,0.75,0.5,0.5,0.25)$ in the domain (\ref{FF2}) and $(c_1,c_2,c_3,c_4,c_5)=(0.5,0.25,0.5,0,0.75)$ in the domain (\ref{FF4}). {\color{black} The green points in Fig.~\ref{fig3a-zoom} is where Mathematica computes different optimal values ($(c_1,c_2,c_3,c_4,c_5)=(0.125,0.625,0.5,0.625,0.125)$ in the domain (\ref{FF2}) and  $(c_1,c_2,c_3,c_4,c_5)=(0.625,0.125,0.5,0.125,0.625)$ in the domain (\ref{FF4})), however, the value of the functional $F_R$ for these points roughly coincides (up to error of $\pm 0.004$) with what we get if consider $(c_1,c_2,c_3,c_4,c_5)=(0,0.75,0.5,0.5,0.25)$ and $(c_1,c_2,c_3,c_4,c_5)=(0.5,0.25,0.5,0,0.75)$ instead.}






\vskip0.2cm

\noindent For the points $(k_1,k_2)$ where $(F,R)=(1,2)$ {\color{black}(i.e. below the threshold of Fig.~\ref{fig3a-zoom})}, the optimal values of $c_i$ are $(c_1,c_2,c_3,c_4,c_5)=(0.5,0.5,0,0.5,0.5)$ (Mathematica computes these values with a $\pm0.02$ error) regardless of whether (\ref{FF2}) or (\ref{FF4}) is used. 

\vskip0.2cm


\begin{figure}[h]
\centering
\includegraphics[scale=0.6]{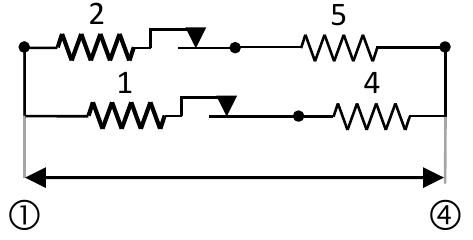}
\caption{The model of Fig.~\ref{fig1} with spring 3 removed (i.e. with $c_3=0$). The bold springs denote the springs that will deform plastically as the displacement-controlled loading (double-sided arrow) increases.} \label{figMG}
\end{figure}

\begin{figure}[h]
\centering
\includegraphics[scale=0.6]{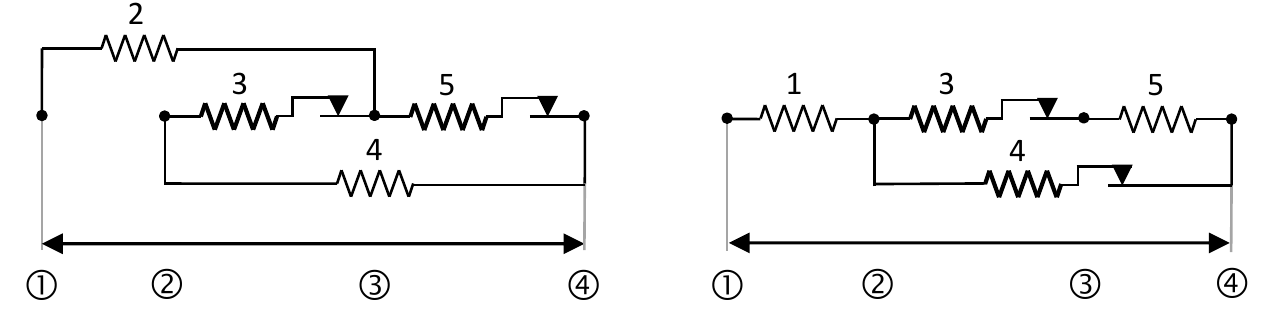}
\caption{The models of Fig.~\ref{fig1} with spring 1 (left) and spring 2 (right) removed (i.e. with $c_1=0$ and $c_2=0$ accordingly). The meaning of bold springs is the same as in Fig.~\ref{figMG}.} \label{figMR}
\end{figure}

\begin{proposition}\label{prop32}  In order to obtain the maximal interplay between strength and resistance of the network of Fig.~\ref{fig1} (functional $M_R$) keeping the cost of the network under a fixed value, only 4 of 5 springs of Fig.~\ref{fig1} are needed for most of the values of $(k_1,k_2)$ except for a few exceptional values of $(k_1,k_2)$ shown in green in Figs.~\ref{fig3a} and \ref{fig3a-zoom} where all 5 springs are needed. In the case where 4 springs are needed, the 4 springs should be connected as in Fig.~\ref{figMG} or as in Fig.~\ref{figMR} according to whether the parameters $(k_1,k_2)$ are above of a certain threshold (shown in Fig.~\ref{fig3a-zoom}) or below. For  parameters  $(k_1,k_2)$ above the threshold and on the threshold, the equal maximums of $F_R$ are attained at 
$(c_1,c_2,c_3,c_4,c_5)=(0,0.75,0.5,0.5,0.25)$ (domain (\ref{FF2})) and  $(c_1,c_2,c_3,c_4,c_5)=(0.25,0.5,0.5,0.75,0)$ (domain (\ref{FF4})) except for  special values of $(k_1,k_2)$ near the threshold (indicated in Fig.~\ref{fig3a-zoom} in green) where equal maximums of $F_R$ are attained at $(c_1,c_2,c_3,c_4,c_5)=(0.125,0.625,0.5,0.625,0.125)$ (domain (\ref{FF2})) and $(c_1,c_2,c_3,c_4,c_5)=(0.625,0.125,0.5,0.125,0.625)$ (domain (\ref{FF4})). 
For parameters  $(k_1,k_2)$ below the threshold, the maximum of $F_R$ is attained at $(c_1,c_2,c_3,c_4,c_5)=(0.5,0.5,0,0.5,0.5).$ These values of $c_i$ belong to the boundaries of both domains (\ref{FF2}) and (\ref{FF4}) simultaneously.
\end{proposition}

\noindent Proposition~\ref{prop32} implies that, if electrical resistance of the material is a more critical objective compared to mechanical strength, then the topology of Fig.~\ref{figMG} has to be chosen for the model design; if mechanical strength is valued more over the electrical resistance then the topology of Fig.~\ref{figMR} needs to be used. 

\vskip0.2cm

\noindent From Fig.~\ref{fig3a} we can also conclude  that $F_R$ increases faster along the direction $(0,1)$ (of the $(k_1,k_2)$ coordinate plane) compared to the direction $(1,0).$

\subsection{Maximization of the strength-conductance functional $F_G$}\label{Sec33}

\noindent The settings for optimization of the functional $F_G$ (given by (\ref{functionals1}))  and Wolfram Mathematica code are analogous to those of functional $F_R$, see the beginning of  Section~\ref{SectionMR}.

\vskip0.2cm

\noindent We conclude that for each $(k_1,k_2)$, the maximal value of $F_G$ is attained at  $c_1=c_2=c_4=c_5=0.5,$  $c_3=0$, and the corresponding values of $F$ and $G$ are  $(F,G)=(1,0.5)$. 

\begin{proposition}  In order to obtain the maximal interplay between strength and conductance of the network of Fig.~\ref{fig1} (functional $F_G$) keeping the cost of the network under a fixed value, only 4 of 5 springs of Fig.~\ref{fig1} are needed and the topology of the network should be taken as in   Fig.~\ref{figMG}. The values of $c_i$ at which the maximum of $F_G$ is attained are $(c_1,c_2,c_3,c_4,c_5)=(0.5,0.5,0,0.5,0.5)$ for all $(k_1,k_2).$ These values of $c_i$ belong to the boundaries of both domains (\ref{FF2}) and (\ref{FF4}) simultaneously. \label{prop33}
\end{proposition}

\begin{figure}[H]
\centering
\includegraphics[clip, trim=0.3cm 8.5cm 0.3cm 8.5cm, width=1.00\textwidth]{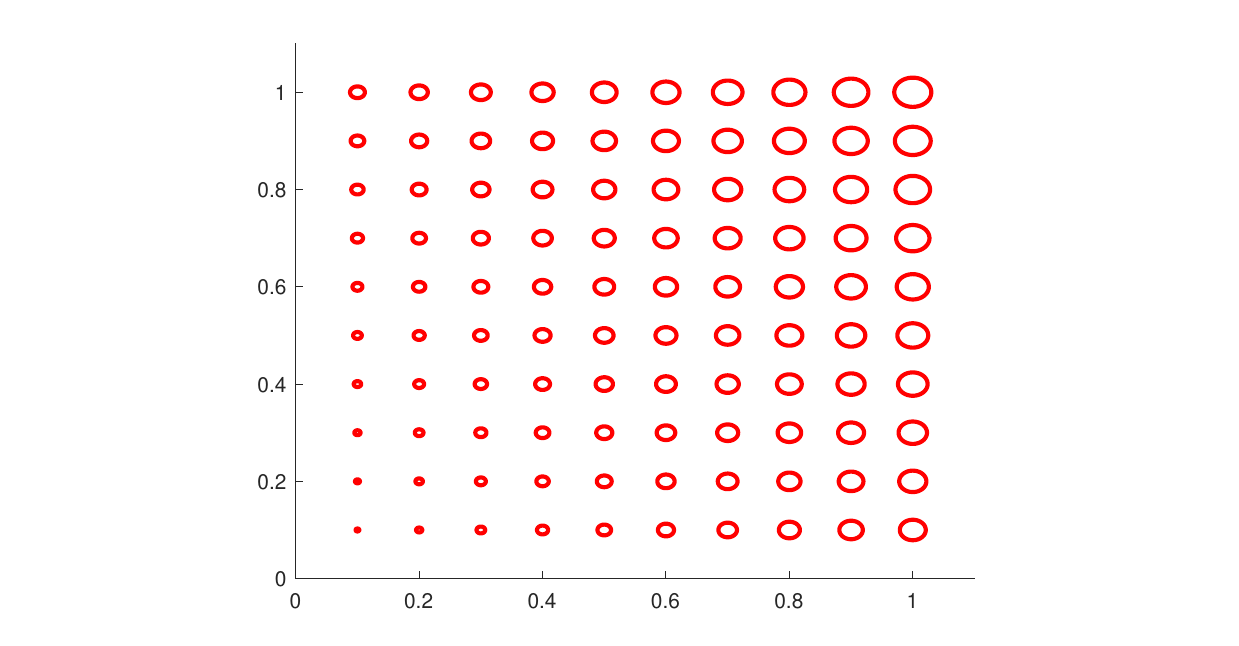}
\caption{Maximization of the functional $F_G$. The horizontal and vertical axis stay for $k_1$ and $k_2$ respectively. The radius of the disc centered at $(k_1,k_2)$ refers to the maximum of $F_R$ for the parameters  $(k_1,k_2)$.  } \label{fig3b}
\end{figure}

\noindent The ratios of how $F_G$ increases across different directions of $(k_1,k_2)$ can be learned from Fig.~\ref{fig3b}. We conclude that $F_G$ increases faster along the direction $(1,0)$ (of the $(k_1,k_2)$ coordinate plane) compared to the direction $(0,1).$

\section{Minimization of the cost of  material under constraint on the multi-functional properties}\label{Sec4}

\noindent Since the optimal values of the functionals $F_R$ and $F_G$ vary in the intervals [0.1,1.2] and [0.1,1.5]  when the cost of the material is constrained (Sections~\ref{SectionMR} and \ref{Sec33}), we will pick a constant inside both of the intervals as a constraint (lower bound) for the functionals $F_R$ and $F_G$ when optimizing $C$. Specifically, we will consider 
\begin{equation}\label{FR-05}
    F_R\ge 0.5
\end{equation}
and 
\begin{equation}\label{FG-05}
    F_G\ge 0.5
\end{equation}
in Sections~\ref{FR05} and \ref{FG05} to follow.

\subsection{Minimization of the cost under constraint on the strength-resistance functional} \label{FR05}
\begin{figure}[H]
\centering
\hskip-4.8cm\includegraphics[scale=0.65]{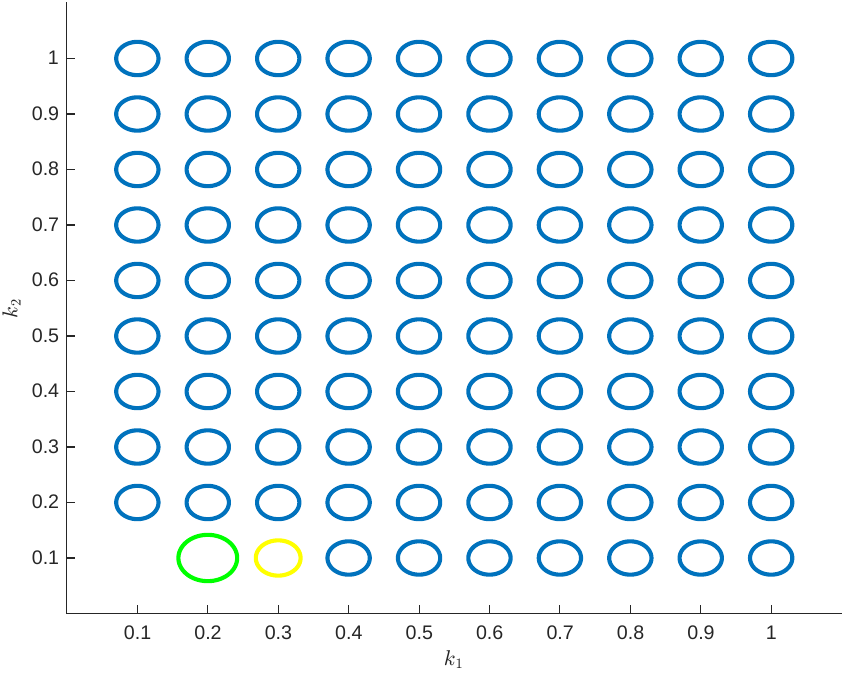}
\vskip-6cm\hskip8.2cm
$
\ \ \ \ \begin{array}{l}
\ \ \ (c_1,\ \ c_2,\ \ c_3,\ \ c_4,\ \ c_5)\\
\\
\mbox{\definecolor{blue1}{rgb}{0 0.4470 0.7410}\textcolor{blue1}{$\blacksquare$}}\ (0.16,0.58,0,0.16,0.58)\\
\mbox{\color{yellow}$\blacksquare$}\ \mbox{other values}\\
\mbox{\color{green}$\blacksquare$}\ \mbox{other values}\\

\end{array}
$
\vskip3.6cm

\caption{Minimization of the functional $C$ over varying $(k_1,k_2)$. 
The radius of the disc centered at $(k_1,k_2)$ refers to the minimum of $C$ for the parameters  $(k_1,k_2)$. 
}\label{fig4anew} 
\end{figure}

\begin{figure}[H]
\centering
\hskip-4.8cm\includegraphics[trim = 0mm 60mm 0mm 70mm,scale=0.65]{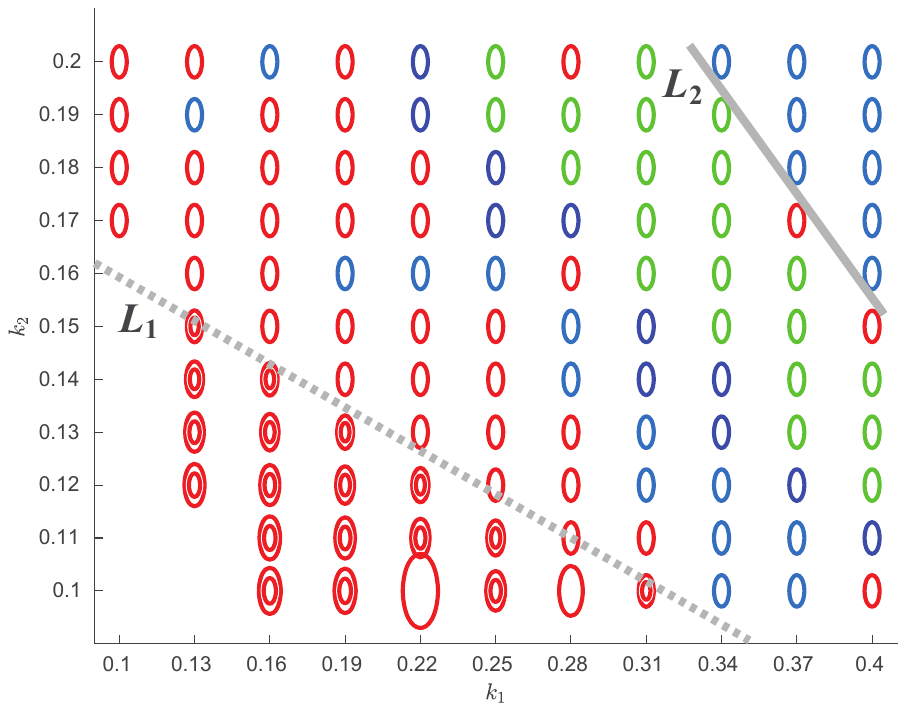}
\vskip-7.15cm\hskip9.2cm
$
\begin{array}{l}
\mbox{\color{white}$\blacksquare$}\  
(c_1,\ \ c_2,\ \ c_3,\ \ c_4,\ \ c_5)
\\
\\
\mbox{\color{blue}$\blacksquare$}\ (0.34,0.40,0,0.34,0.40)\\
\mbox{\color{green}$\blacksquare$}\ (0.51,0.23,0,0.51,0.23)\\
\mbox{\definecolor{nblue}{rgb}{0 0.4470 0.7410}\textcolor{nblue}{$\blacksquare$}}\ (0.59,0.15,0,0.59,0.15)\\
\mbox{\color{red}$\blacksquare$}\ (\text{different values})
\end{array}
$
\vskip3.5cm
\caption{Minimization of the functional $C$. The horizontal and vertical axis stay for $k_1$ and $k_2$ respectively. The radius of the disc centered at $(k_1,k_2)$ refers to the minimum of $C$ for the parameters  $(k_1,k_2)$. For each color we have a different optimal value of $C$. The circles in red represent the case where optimal value of $C$ and $(c_1,c_2,c_3,c_4,c_5)$ are different for each case. The concentric circle corresponds to the cases where $c_3 >0$ and also where optimal value of $(c_1,c_2,c_3,c_4,c_5)$ does not satisfy the feasibility conditions for domain (\ref{FF2}).
}\label{newfig1} 
\end{figure}

\begin{figure}[H]
\centering
\hskip-4.8cm\includegraphics[trim = 0mm 70mm 0mm 80mm,scale=0.65]{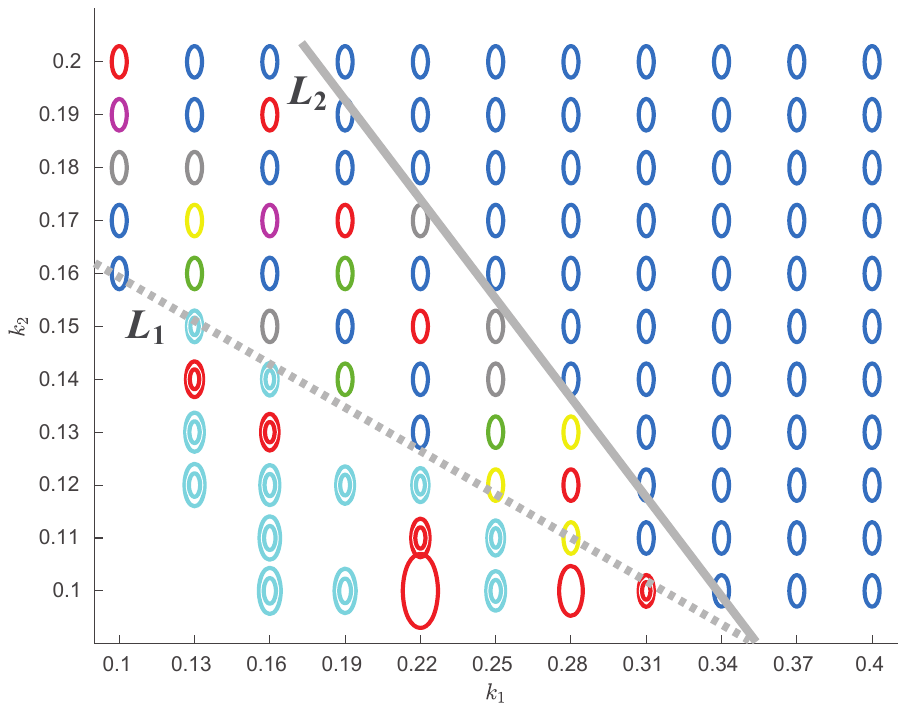}
\vskip-6.15cm\hskip9.2cm
$
\begin{array}{l}
\mbox{\color{white}$\blacksquare$}\  
(c_1,\ \ c_2,\ \ c_3,\ \ c_4,\ \ c_5)
\\
\\
\mbox{\definecolor{cyan}{rgb}{0 1 1}\textcolor{cyan}{$\blacksquare$}}\ (0.75,\mbox{various values})\\
\mbox{\definecolor{ngreen}{rgb}{0.4660 0.6740 0.1880}\textcolor{ngreen}{$\blacksquare$}}\ (0.26,0.48,0,0.26,0.48)\\
\mbox{\definecolor{nblue}{rgb}{0 0.4470 0.7410}\textcolor{nblue}{$\blacksquare$}}\ (0.59,0.15,0,0.59,0.15)\\
\mbox{\color{yellow}$\blacksquare$}\ (0.54,0.20,0,0.54,0.20)\\
\mbox{\definecolor{magenta}{rgb}{1,0,1}\textcolor{magenta}{$\blacksquare$}}\ 
(0.25,0.49,0,0.25,0.49)\\
\mbox{\definecolor{gray}{rgb}{0.5019 0.5019 0.5019}\textcolor{gray}{$\blacksquare$}}\ 
(0.11,0.63,0,0.11,0.63)\\
\mbox{\color{red}$\blacksquare$}\ (\text{various values})\\
\end{array}
$
\vskip1.75cm
\caption{Minimization of the functional $C$. The horizontal and vertical axis stay for $k_1$ and $k_2$ respectively. The radius of the disc centered at $(k_1,k_2)$ refers to the minimum of $C$ for the parameters  $(k_1,k_2)$. For each color we have a different optimal value of $C$. The circles in red represent the case where optimal value of $C$ and $(c_1,c_2,c_3,c_4,c_5)$ are different for each case. The concentric circle corresponds to the cases where $c_3 >0$ and also where optimal value of $(c_1,c_2,c_3,c_4,c_5)$ does not satisfy the feasibility conditions for domain (\ref{FF4}).
}\label{newfig2} 
\end{figure}

 For each of the two domains (\ref{FF2}) and (\ref{FF4}), and any $(k_1,k_2)\not=(0.1,0.1)$, the minimal value of  $C$ with the constraints (\ref{optF}) and (\ref{FR-05}) is always attained on the boundary. 
 As Fig.~\ref{fig4anew} illustrates, optimal values of $c_i$ and $C$ don't change for almost the entire range of parameters $(k_1,k_2)$ with some changes taking place closer to $(k_1,k_2)=(0.1,0.1).$ 
To get insight into these changes we zoom the region $(k_1,k_2)\in[0.1,0.4]\times[0.1,0.2]$ in Figs.~\ref{newfig1} and \ref{newfig2} (corresponding to domains (\ref{FF2}) and (\ref{FF4}) accordingly).

\vskip0.2cm

\noindent From Figures~\ref{newfig1} and \ref{newfig2} we observe that the optimal values $(c_1,...,c_5)$ and the cost value $C$ stay constant in the upper right part of the $(k_1,k_2)$ plane, specifically
\begin{eqnarray}
 && (c_1,...,c_5)=(0.594,0.156,0,0.594, 0.156)\quad\mbox{in domain (\ref{FF2})},\\
 && (c_1,...,c_5)=(0.162,0.588,0,0.162, 0.588)\quad\mbox{in domain (\ref{FF4})}
\end{eqnarray}
(which differ by $\pm0.006$ that we relate to computational error), 
until $(k_1,k_2)$ crosses a straight line (threshold) of an approximate slope of $-0.5$ (solid gray lines in Figs.~\ref{newfig1} and \ref{newfig2}). Moving further towards $(k_1,k_2)=(0.1,0.1)$, we observe another threshold (dotted gray lines in Figs.~\ref{newfig1} and \ref{newfig2} close to $-0.5$ slope-wise too) when $c_3$ switches from zero to positive. In between the solid and dotted lines, we discovered several collections of equal $c=(c_1,...,c_5)$ indicated in Figs.~\ref{newfig1} and \ref{newfig2} in different colors. 
As follows from Figs.~\ref{fig4anew}, \ref{newfig1}, \ref{newfig2}, at least one $c_i$ vanishes in each of the optimal values $(c_1,...,c_5)$. 
In particular, in the case of domain (\ref{FF2}), non-vanishing $c_3$ implies that either $c_1=0$ or $c_5=0$ (i.e. the topology of the network takes the form of Fig.~\ref{figMR}\,(left)), while in the case of domain (\ref{FF4}), non-vanishing $c_3$ implies that either $c_2=0$ or $c_4=0$ (i.e. the topology of the network takes the form of Fig.~\ref{figMR}\,(right)).

\vskip0.2cm

\noindent We discover that the minimal value of $C=1.5$ for the values of $(k_1,k_2)$ above the dotted threshold and $C>1.5$ below the dotted threshold, i.e. whether $C=1.5$ or $C>1.5$ turns out to be strictly related to whether $c_3=0$ or $c_3>0$ with 4 exceptions listed in Tables~\ref{table1} and \ref{table2}. 

\begin{table}[!h]
\centering
\begin{tabular}{ |c|c|c|c|c| } 
 \hline
 $(k_1,k_2)$ & $(c_1,c_2,c_3,c_4,c_5)$ & $F$ & $R$ & $C$ \\ 
 \hline
  $(0.22,0.1)$ & $(0.88,0.88,0,0.88,0.88)$ & 1.75 & 1.14 & 3.51 \\ 
  $(0.28,0.1)$ & $(0.46,0.72,0,0.46,0.72)$ & 1.18 & 1.69 & 2.36 \\ 
 \hline
\end{tabular}
\vspace{0.5em} 
\centeredcaption{The 2 cases where $C>1.5$ doesn't imply $c_3>0$ in domain (\ref{FF2}).\label{table1}}
\end{table}

\begin{table}[!h]
\centering
\begin{tabular}{ |c|c|c|c|c| } 
 \hline
 $(k_1,k_2)$ & $(c_1,c_2,c_3,c_4,c_5)$ & $F$ & $R$ & $C$ \\ 
 \hline
   $(0.22,0.1)$ & $(1.15,0.6,0,1.15,0.6)$ & 1.75 & 1.14 & 3.51 \\ 
  $(0.28,0.1)$ & $(0.46,0.72,0,0.46,0.72)$ & 1.18 & 1.69 & 2.36 \\ 
 \hline
\end{tabular}
\vspace{0.5em}
\centeredcaption{The 2 cases where $C>1.5$ doesn't imply $c_3>0$ in domain (\ref{FF4}).\label{table2}}
\end{table}


\noindent Notice, the values of $(c_1,c_2,c_3,c_4,c_5)$ in Tables~\ref{table1} and \ref{table2} is what we obtain with Wolfram Mathematica. However, one can immediately see that each of the $(c_1,c_2,c_3,c_4,c_5)$ in Tables~\ref{table1} and \ref{table2} satisfy both domains (\ref{FF2}) and (\ref{FF4}) at the same time.

\vskip0.2cm

\noindent The optimal value of the functional $F_R=0.5$ below the dotted threshold and $F_R>0.5$ above the dotted threshold. The corresponding value of response force is always $F=0.75$ except for the cases listed in Tables~\ref{table1} and \ref{table2}.  

\begin{proposition}\label{prop41}  In order to minimize the manufacturing cost of the network of Fig.~\ref{fig1} keeping both the strength and strength-resistance above fixed values, only 4 of 5 springs of Fig.~\ref{fig1} are needed (connected as in   Fig.~\ref{figMG} above the dotted threshold $L_1$ of Figs.~\ref{newfig1}-\ref{newfig2} and connected as in Fig.~\ref{figMR} below $L_1$). 
The optimal values are approximately 
$
(c_1,...,c_5)=(0.59,0.16,0,0.59,0.16)
$
above the solid threshold $L_2$ and are grouped in patterns shown in Figs.~\ref{newfig1}-\ref{newfig2} below $L_2$. The minimal value of $C=1.5$ above $L_1$ and $C>1.5$ below $L_1$ except for a few exceptions that occur at $k_2=0.1$ (and that are listed in Tables~\ref{table1} and \ref{table2}). The values $F_R>0.5$ and $F=0.75$ above $L_1$ suggest the cost $C=1.5$ above $L_1$ is determined solely by constraint (\ref{optF}) and constraint (\ref{FR-05}) is redundant above $L_1.$ On the other hand $F_R=0.5$ and $F=0.75$ below $L_1$ suggests that both constraints (\ref{optF}) and (\ref{FR-05}) come to play below $L_1$. The slopes of both thresholds $L_1$ and $L_2$ are approximately $-1/2.$ For each of the domains (\ref{FF2}) and (\ref{FF4}), the maximal value of the minimum of $C$ (the biggest oval in Figs.~\ref{newfig1} and \ref{newfig2}) is unique across all parameters $(k_1,k_2)$ and attained at approximately $(k_{1,*},k_{2,*})=(0.22,0.1).$ The optimal values $(c_1,...,c_5)$ are non-unique, but the threshold $L_1$ doesn't depend on the domain. {\color{black}While for $(k_1,k_2)$ above $L_1$ the optimal values of $(c_1,...,c_5)$ on (\ref{FF2}) coincide with the optimal values of $(c_1,...,c_5)$ on (\ref{FF4}) (i.e. common optimal values $(c_1,...,c_5)$ belong to the boundaries of both domains simultaneously), for values of $(k_1,k_2)$ below $L_1$, the optimal values $(c_1,...,c_5)$ on (\ref{FF2}) are different from the optimal values $(c_1,...,c_5)$ on  (\ref{FF4}).}
\end{proposition}

\vskip0.2cm

\subsection{ Minimization of the cost under constraint on the strength-conductance functional }\label{FG05}
\begin{figure}[H]
\centering
\hskip-4.8cm\includegraphics[scale=0.65]{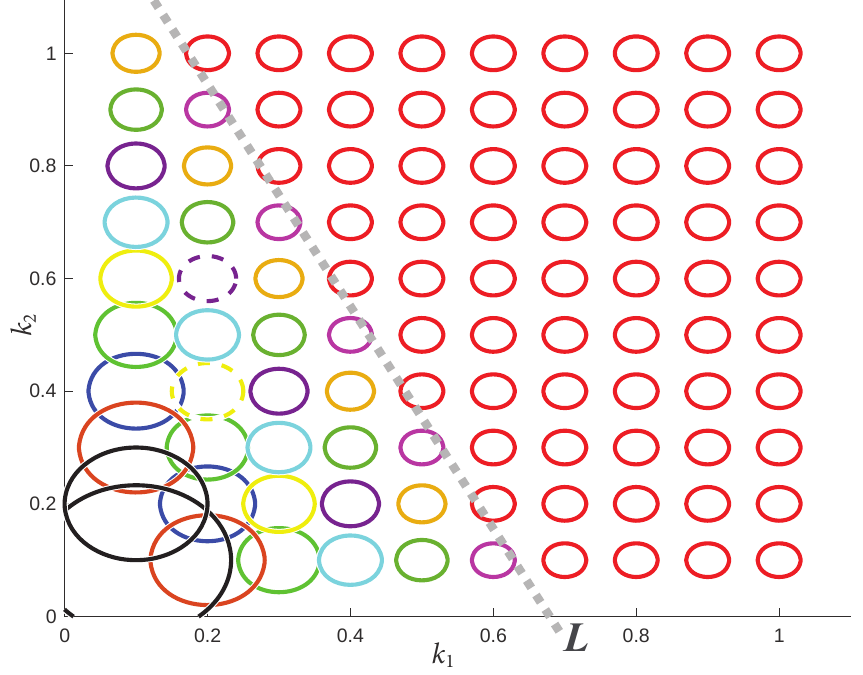}
\vskip-7.15cm\hskip9.2cm
$
\begin{array}{l}
\mbox{\color{white}$\blacksquare$}\  
(c_1,\ \ c_2,\ \ c_3,\ \ c_4,\ \ c_5)
\\
\\
\mbox{\color{blue}$\blacksquare$}\ (0.83,0.83,0,0.83,0.83)\\
\mbox{\color{green}$\blacksquare$}\ (0.71,0.71,0,0.71,0.71)\\
\mbox{\color{yellow}$\blacksquare$}\ (0.62,0.62,0,0.62,0.62)\\
\mbox{\color{yellow}$\square$}\ (0.59,0.65,0,0.59,0.65)\\
\mbox{\definecolor{cyan}{rgb}{0,1,1}\textcolor{cyan}{$\blacksquare$}}\ (0.55,0.55,0,0.55,0.55)\\
\mbox{\definecolor{purple}{rgb}{0.4940, 0.1840, 0.5560}\textcolor{purple}{$\blacksquare$}}\ (0.5,0.5,0,0.5,0.5)\\
\mbox{\definecolor{purple}{rgb}{0.4940, 0.1840, 0.5560}\textcolor{purple}{$\square$}}\ (0.58,0.41,0,0.58,0.41)\\
\mbox{\definecolor{dullgreen}{rgb}{0.4660 ,0.6740, 0.1880}\textcolor{dullgreen}{$\blacksquare$}}\ (0.45,0.45,0,0.45,0.45)\\
\mbox{\definecolor{mustard}{rgb}{0.9290 ,0.6940, 0.1250}\textcolor{mustard}{$\blacksquare$}}\ (0.41,0.41,0,0.41,0.41)\\
\mbox{\definecolor{magenta}{rgb}{1,0,1}\textcolor{magenta}{$\blacksquare$}}\ 
(0.38,0.38,0,0.38,0.38)\\
\mbox{\definecolor{orange}{rgb}{0.8500 ,0.3250, 0.0980}\textcolor{orange}{$\blacksquare$}}\ (1,1,0,1,1)
\end{array}
$
\vskip1.3cm
\caption{Minimization of the cost $C$ the strength-conductance functional $F_G$. The radius of the disc centered at $(k_1,k_2)$ refers to the minimum of $C$ for the parameters  $(k_1,k_2)$. The
optimal value $C=1.5$ or $C>1.5$ according to whether $(k_1,k_2)$ is above (red circles) or below (non-red circles) the solid threshold. The values of $c$ for different colors are listed in the legend. These values of $c$ appear to not depend on whether domain (\ref{FF2}) and (\ref{FF4}) is considered except for the dashed circles: (i) for $(k_1,k_2)$ corresponding to the dashed purple circle, the value of $c$ coincides with $c$ values for other purple circles only for the case of domain (\ref{FF2}) (the value of $c$ for the domain (\ref{FF4} is noted in the legend with a hollow rectangle), (ii) for $(k_1,k_2)$ corresponding to the dashed yellow circle, the value of $c$ coincides with $c$ values for other yellow circles only for the case of domain (\ref{FF4}) (the value of $c$ for the domain (\ref{FF2} is noted in the legend with a hollow rectangle).    We do not see any distinguishable patterns of equal $c$ for the values of $(k_1,k_2)$ that correspond to the red circles. The circles in black just represent all the remaining cases. }\label{fig4a} 
\end{figure}

 For each of the two domains (\ref{FF2}) and (\ref{FF4}), and any $(k_1,k_2)\not=(0.1,0.1)$, the minimal value of  $C$ with the constraints (\ref{optF}) and (\ref{FG-05}) is always attained on the boundary.  As Fig.~\ref{fig4a} illustrates, optimal value $C=1.5$ or $C>1.5$ according to whether $(k_1,k_2)$ is above or below the dotted threshold (of approximate slope $-2$).  For all $(k_1,k_2)$, the optimal values of $c=(c_1,...,c_5)$ satisfy either 
 
 (i) $c_1=c_2=c_4=c_5$, $c_3=0$ or
 
 (ii) $c_1=c_4$, $c_2=c_5,$ $c_3=0$.

\noindent In particular, the optimal value of $c_3$ is always 0.
  For the values $(k_1,k_2)$ below the dotted threshold, we observe lines (with the slope roughly equal to the slope of the dotted threshold) of equal optimal values of $c$, which can be seen in different colors in Fig.~\ref{fig4a}. At the same time, we do not see any distinguishable patterns of equal optimal values of $c$ when $(k_1,k_2)$ is above the solid threshold.

\vskip0.2cm

\noindent For all the cases where the minimal cost is $C=1.5$ (red circles in the graph), we get $F=0.75$, $G=0.375$ and $F_G>0.5$. For all the cases where $C>1.5$ (non-red circles) we get $F>0.75,$ $G>0.375$, but, interestingly, $F_G=0.5$. 



\begin{proposition}  In order to minimize the manufacturing cost of the network of Fig.~\ref{fig1} keeping both the strength and strength-conductance above fixed values (constraints (\ref{optF}) and (\ref{FG-05})), only 4 of 5 springs of Fig.~\ref{fig1} are needed connected as in Fig.~\ref{figMG}. 
For each of the domains (\ref{FF2}) and (\ref{FF4}), the minimal value of $C=1.5$ above the dotted threshold $L$ (of approximate slope $-2$) of Fig.~\ref{fig4a} and $C>1.5$ below $L$. The value $F=0.75$ above $L$ suggests  that constraint (\ref{optF}) doesn't allow the cost to go below $C=1.5$ for $(k_1,k_2)$ above $L$ and constraint (\ref{FG-05}) doesn't come to play above $L.$ On the other hand, $F_G=0.5$ below $L$ suggests that 
constraint (\ref{optF}) is redundant below $L$. For each $(k_1,k_2)$, {\color{black}the optimal values $(c_1,c_2,c_3,c_4,c_5)$ belong to the boundaries of both domains (\ref{FF2}) and (\ref{FF4}) simultaneously} (and $c_3$ is always 0). The other optimal parameters $c_1,c_2,c_4,c_5$ form lines (of approximate slope $-2$) of equal values below the threshold, as shown in Fig.~\ref{fig4a}. The minimal value of $C$ increases as $(k_1,k_2)$ approaches $(k_1,k_2)=(0.1,0.1)$. \label{prop42}
\end{proposition}



\section{Conclusions}
In this paper we considered a prototypic model of 5 elastoplastic conducting springs on 4 nodes and used formulas from the sweeping process theory to determine  parameters of the springs that solve the following optimization problems: 
\begin{itemize}
\item[(A)] maximize the combined strength  and electric resistance of material (functional $F_R$) under constraint on the maximal fabrication cost $C$ and strength $F$; 
\item[(B)] maximize the combined strength  and electric conductance of the material (functional $F_G$) under constraint on the maximal values of $C$ and strength $F$;
\item[(C)] minimize the fabrication cost $C$ assuming that neither $F_R$, nor $F$,  go below given constraints; 
\item[(D)] minimize the fabrication cost $C$ assuming that neither $F_G$, nor $F$, go below given constraints. 
\end{itemize}

\noindent The conclusions obtained with Wolfram Mathematica NMaximize command are summarized in Propositions~\ref{prop32}, \ref{prop33}, \ref{prop41}, \ref{prop42} and motivate several tasks towards rigorous mathematical proofs that we group below according to the above stated Problems A, B, C, and D. We recall that we refer to an optimal value $(c_1,...,c_5)$ as unique, if $(c_1,...,c_5)$ is unique in one of the domains (\ref{FF2}) and (\ref{FF4}), see Remark~\ref{uniquenessremark}.


\vskip0.2cm

\noindent A1: Prove uniqueness of the optimal value $(c_1,...,c_5)$ for each value of $(k_1,k_2).$

\vskip0.2cm

\noindent A2: Prove that the space of parameters $(k_1,k_2)$ splits into 3 groups so that the optimal values $(c_1,...,c_5)$ in each group coincide. 

\vskip0.2cm

\noindent 
A3: Examine whether or not one of the 3 groups from A1 always forms a straight line that separates the other 2 groups. At the very least, prove that the above-mentioned 2 groups are always separated by a straight line $L$. Find the equation of $L$, which is expected to be independent of the chosen domain of parameters $(c_1,...,c_5)$ (domains (\ref{FF2}) and (\ref{FF4}) in this paper).

\vskip0.2cm

\noindent A4: Prove that the optimal values $(c_1,...,c_5)$ below $L$ belong to the boundaries of both domains (\ref{FF2}) and (\ref{FF4}) simultaneously. Prove that the optimal values $(c_1,...,c_5)$ above $L$ and on $L$ do not belong to the boundaries of (\ref{FF2}) and (\ref{FF4}) simultaneously. 

\vskip0.2cm

\noindent 
B1: Prove that optimal values $(c_1,...,c_5)$ are unique, independent of parameters $(k_1,k_2),$ and belong to the boundaries of both domains (\ref{FF2}) and (\ref{FF4}) simultaneously.

\vskip0.2cm

\noindent 
C1: Prove the existence of a linear threshold $L_2$ in the space of parameters $(k_1,k_2)$ on so that the optimal values $(c_1,...,c_5)$ above $L_2$ coincide. Find the equation of $L_2$, which is expected to depend on the chosen domain of parameters $(c_1,...,c_5)$.

\vskip0.2cm

\noindent 
C2: For the values of $(k_1,k_2)$ below $L_2$, prove the existence of lines (presumably parallel to $L_2$) so that the optimal values of $(c_1,...,c_5)$ are constant alone these lines.

\vskip0.2cm

\noindent 
C3: Prove the existence of another linear threshold $L_1$ in the space of parameters $(k_1,k_2)$ (below $L_2$) such that $C$ keeps a constant value above $L_1$ and $C$ is larger than this constant value below $L_1.$ Find the equation of $L_1$, which (in contrast with $L_2$) is expected to be independent on the chosen domain of parameters $(c_1,...,c_5$ and expected to be parallel to $L_2$. 

\vskip0.2cm

\noindent 
C4: Prove that crossing the threshold $L_1$ corresponds to $F_R$ reaching its chosen bound {\color{black}(in this paper, crossing $L_1$ from above to below, $F_R$ reaches the lower bound of $0.5$)}.  

\vskip0.2cm

\noindent 
C5: Prove that $F$ equals its chosen bound (0.75 in this paper) for all $(k_1,k_2)$ except for isolated special values. Determine the nature of these special values.

\vskip0.2cm

\noindent 
C6: Prove that, across all values of $(k_1,k_2)\in[k_{1,min},k_{1,max}]\times[k_{2,min},k_{2,max}]$, the {\color{black}maximum of the minimal values}  of $C$ is attained at some $(k_{1,*},k_{2,*})\not=(k_{1,min},k_{2,min}).$ Moreover, such a value of $(k_{1,*},k_{2,*})$ is unique.

\vskip0.2cm

\noindent C7: Prove the existence of a linear threshold such that the optimal values $(c_1,...,c_5)$ are unique above the threshold and non-unique below the threshold.



\vskip0.2cm

\noindent {\color{black}C8: Prove that the optimal values $(c_1,...,c_5)$ above $L_1$ belong to the boundaries of both domains (\ref{FF2}) and (\ref{FF4}) simultaneously. Prove that the optimal values $(c_1,...,c_5)$ below $L_1$ do not belong to the boundaries of (\ref{FF2}) and (\ref{FF4}) simultaneously. }

\vskip0.2cm

\noindent 
D1: Prove the existence of a linear threshold $L$ in the space of parameters $(k_1,k_2)$ such that $C$ keeps a constant value above $L$ and $C$ is larger than this constant value below $L.$ Find the equation of $L$, which is expected to be independent of the chosen domain of parameters $(c_1,...,c_5)$.

\vskip0.2cm

\noindent 
D2: Prove that crossing the threshold $L$ corresponds to $F$ switching from $F>const$ to $F=const$ ($const=0.75$ in this paper) and switching from $F_G=const$ to $F_G>const$ ($const=0.5$), if moving $(k_1,k_2)$ from under $L$ to above $L$.

\vskip0.2cm

\noindent 
D3: For the values of $(k_1,k_2)$ below $L$, prove the existence of lines along where the optimal values of $(c_1,...,c_5)$ are constant.

\vskip0.2cm

\noindent 
D4: Prove uniqueness of the optimal values $(c_1,...,c_5)$ for every $(k_1,k_2)$ potentially disproving the existence of those special pairs $(k_1,k_2)$ (such as $(k_1,k_2)=(0.2,0.6)$ and $(k_1,k_2)=(0.2,0.4)$ in Fig.~\ref{fig4a}) where nonuniqueness is detected by Wolfram Mathematica (or explain the nature of $(k_1,k_2)$ where nonuniquness of the optimal values takes places).

\vskip0.2cm

\noindent Simulations of this paper are executed under assumption that the plastic mode is terminally reached by a maximally possible number of springs. Dropping this assumption might be of interest in applications outside of materials science. This would require extension of our numerical analysis from domains (\ref{FF2}) and (\ref{FF4}) to all positive values of $c_1,...,c_5$, that might stimulate further development of the theory of \cite{SICON} and further derivation of the corresponding terminal forces  (\ref{FF1}) and (\ref{FF4}) (that are not yet available for all positive values of $c_1,...,c_5$).

\vskip0.2cm

\noindent We hope that simulation results reported and theoretical questions motivated by these simulation results will attract interest of future researchers.

\section{Appendix: Conductance of a Wheatstone Bridge} \label{appA}

This section follows the lines of  \cite{resistance} to calculate the equivalent conductance of the circuit of Fig.~\ref{fig1} where spring $i$ is viewed as resistance $R_i$. We assume that node 1 is connected to the positive terminal of the battery of potential $\epsilon$ and node 4 is connected to the negative terminal of zero potential. As per the junction law, the sum of currents into node 2 and node 3 is $0$. 
$$\frac{\epsilon-V_{2}}{R_{1}}+\frac{V_{3}-V_{2}}{R_{3}}+\frac{0-V_{2}}{R_{4}}=0$$
$$\frac{\epsilon-V_{3}}{R_{2}}+\frac{V_{2}-V_{3}}{R_{3}}+\frac{0-V_{3}}{R_{5}}=0$$
We consider $v_{2}=\frac{V_{2}}{\epsilon}$ and $v_{3}=\frac{V_{3}}{\epsilon}$, so we get , 

$$\left(\frac{1}{R_{1}}+\frac{1}{R_{3}}+\frac{1}{R_{4}}\right)v_{2}-\frac{v_{3}}{R_{3}}=\frac{1}{R_{1}}$$

$$\left(\frac{1}{R_{2}}+\frac{1}{R_{3}}+\frac{1}{R_{5}}\right)v_{3}-\frac{v_{2}}{R_{3}}=\frac{1}{R_{2}}$$
$v_{2}$ and $v_{3}$ can be expressed as :
$$\frac{1}{v_{2}}=\frac{R_{3}R_{4}R_{5}+R_{1}(R_{3}+R_{4})R_{5}}{R_{4}((R_{1}+R_{3})R_{5}+R_{2}(R_{3}+R_{5}))}+\frac{R_{2}R_{4}(R_{3}+R_{5})+R_{1}R_{2}(R_{3}+R_{4}+R_{5})}{R_{4}((R_{1}+R_{3})R_{5}+R_{2}(R_{3}+R_{5}))}$$
$$\frac{1}{v_{3}}=\frac{R_{3}R_{4}R_{5}+R_{1}(R_{3}+R_{4})R_{5}}{R_{5}((R_{2}+R_{3})R_{4}+R_{1}(R_{3}+R_{4}))}+\frac{R_{2}R_{4}(R_{3}+R_{5})+R_{1}R_{2}(R_{3}+R_{4}+R_{5})}{R_{5}((R_{2}+R_{3})R_{4}+R_{1}(R_{3}+R_{4}))}$$
The equivalent conductance of the circuit denoted by $\sigma_{eq}$ is the current flowing from the positive to the negative terminal of the battery, given by
$$\sigma_{eq}=\frac{\frac{V_{3}-0}{R_{5}}+\frac{V_{2}-0}{R_{4}}}{\epsilon}=\frac{v_{3}}{R_{5}}+\frac{v_{2}}{R_{4}}$$\\
The equivalent resistance of the system denoted by $R_{eq}$ is given by  
$$R_{eq}=\frac{1}{\sigma_{eq}}=\frac{R_{3}R_{4}R_{5}+R_{1}(R_{3}+R_{4})R_{5}}{R_{3}(R_{4}+R_{5})+(R_{1}+R_{2})(R_{3}+R_{4}+R_{5})}+$$
$$\hspace{1.8cm}\frac{R_{2}R_{4}(R_{3}+R_{5})+R_{1}R_{2}(R_{3}+R_{4}+R_{5})}{R_{3}(R_{4}+R_{5})+(R_{1}+R_{2})(R_{3}+R_{4}+R_{5})}.$$

\section{Appendix: Sweeping process framework for finite-time stability of  elastoplastic systems with application to the benchmark example} \label{appB}

\subsection{The outline of the general framework}
 According to Moreau \cite{Moreau} a network $(D,A,C,R,l(t))$ of $m$ elastoplastic springs on $n$ 1-dimensional nodes with one displacement-controlled loading (we refer the reader to \cite{Gud} for the case of $n$-dimensional nodes) is fully defined by an $m\times n$ kinematic matrix $D$ of the topology of the network, $m\times m$ matrix of stiffnesses (Hooke's coefficients) $A={\rm diag}(a_1,...,a_m)$,  an $m$-dimensional hyperrectangle $C=\prod_{j=1}^m [c_j^-,c_j^+]$ of the achievable stresses of springs (beyond which plastic deformation begins), a vector $R\in\mathbb{R}^m$ of the location of the displacement-controlled loading, and a scalar function $l(t)$ that defines the magnitude of the displacement-controlled loading. Based on the parameters $(D,A,C,R,l(t))$ the space $\mathbb{R}^m$ is decomposed into a direct product of two supspaces $U,V\subset\mathbb{R}^m$ which can be used to formulate a dynamical system that governs the stress vector $s(t)$ of the network. A concise step-by-step guide for computation of all the quantities listed above is available in \cite{PhysicaD}.

\vskip0.2cm

\noindent 
{\bf Sweeping process.} To define a dynamical system that governs the stress-vector of $(D,A,C,R,l(t))$, let 
$$
  N_C^A(y)=\left\{\begin{array}{ll}\left\{\xi\in V:\left<\xi,A(c-y)\right>\le 0,\  c\in C\right\},&\quad {\rm if}\ y\in C,\\
   \emptyset,& \quad {\rm if}\ y\not\in C,
\end{array}\right.
$$
denote the normal cone to a set $C\subset V$ at point $y$. The stress-vector $s(t)$ of 
$(D,A,C,R,l(t))$
is then
\begin{equation}\label{s(t)form}
    s(t)=Ay(t)-l(t)Av,
\end{equation}
where $y(t)$ is the solutions of the differential inclusion (called {\it sweeping process})
\begin{equation}\label{sp}
\begin{array}{c}   -y
   '\in N_{{C}+l(t)v}^A(y), \quad y\in {V},
   \end{array}
\end{equation}
and $v$ is a suitably constructed vector of $V$  (see \cite{PhysicaD,SICON} for the formula of $c(t)=l(t)v$).

\vskip0.2cm

\noindent {\bf Admissibility condition.} To understand how the sweeping process approach helps to predict the asymptotic behavior of $s(t),$ we recall that the polyhedron $C$ in (\ref{sp}) computes as 
$$
{C}=\bigcap\limits_{j\in\overline{1,m}}\left(L(-1,j)\cap L(1,j)\right), 
$$
where  
$$
  L(\alpha,j)=\{x\in V:\left<\alpha n_j,Ax\right>\le \alpha c^\alpha_j\},\qquad (\alpha,j)\in\{-1,1\}\times\overline{1,m}, 
$$
and $n_j$ are suitable vectors of $V$, whose computational formulas (in terms of the mechanical parameters of the network) are also available in \cite{PhysicaD}.
Let $e_j$ be the  vector whose $j$-th component equal $1$ and all other components equal 0.
According to \cite{SICON}, if relation 
\begin{equation}\label{blueinclusion}
\left(\hskip-0.1cm\begin{array}{c}
1  \\
0_{m-n+1}\end{array}\hskip-0.1cm\right) \in {\rm cone}\left(
\left(\hskip-0.1cm\begin{array}{c}
  R^T \\
  (D^\perp)^T\end{array}\hskip-0.1cm\right) \left\{\alpha e_j:(\alpha,j)\in I_0\right\}\right), 
\end{equation}
holds for an {\it irreducible} $I_0\subset\{-1,1\}\times\overline{1,m}$ (i.e. (\ref{blueinclusion}) fails for any $\tilde I_0\subset I_0$), 
where $D^\perp$ is a certain orthogonal to $D$ matrix (see \cite{PhysicaD} for the formula), then  
\begin{equation}\label{Fdef}
\begin{array}{c}     Y(t)=F+l(t)v, \quad{\rm where}\quad
 F= \left(\bigcap\limits_{(\alpha,j)\in I_0} \partial L(\alpha,j)\right)\bigcap {C}, \end{array}
\end{equation}
where $\partial L(\alpha,j)$ is the boundary of $L(\alpha,j)$,
 is a candidate finite-time attractor for sweeping process (\ref{sp}). If $I_0$ satisfies (\ref{blueinclusion}), we say that $I_0$ is {\it admissible}. Therefore, if (\ref{blueinclusion}) holds, then combining (\ref{s(t)form}) and (\ref{Fdef}),
 \begin{equation}\label{SAF}
         S=AF
  \end{equation}
 is a candidate attractor for the stress-vector $s(t)$ of elastoplastic system \\$(D,A,C,R,l(t))$.

\vskip0.2cm

\noindent {\bf Feasibility condition.} To determine, which of the candidate attractors is actually attained, the facet $F$ has to verify a so-called {\it feasibility condition}, i.e. the condition $F\subset C.$ To verify the latter condition computationally, the paper \cite{SICON} determines the vertices of $F$:
\begin{equation}\label{y*iformula}
   y_{*,i}={V}_{basis}\left(
   \left(\left\{e_j,(\alpha,j)\in I_0\cup I_i\right\}\right)^T      
    A V_{basis}\right)^{-1}\left(\left\{c_j^\alpha,(\alpha,j)\in I_0\cup I_i\right\}\right)^T,
\end{equation}
where the columns of $V_{basis}$ are basis vectors of $V$ and $I_i\subset\{-1,1\}\times\overline{1,m}$ is such that $|I_0\cup I_i|=\dim V.$ When $|I_0|=\dim V$, $F$ is just a vertex, i.e. no computation of additional $I_i$ is needed and formula (\ref{y*iformula}) is used with just $M=0$. When $|I_0|<\dim V$, the dimension of $F$ is greater than $0$, i.e. $F$ is a polyhedron whose number of vertices is denoted by $M.$ In terms of the vertices $y_{*,i}$, the {\it feasibility condition} $F\subset C$ is found in \cite{SICON} as follows:
\begin{equation}\label{y*ifeasible}
\begin{array}{ll}
|I_0|<d: &\ \    c_j^-<\left<e_j,Ay_{*,i}\right>< c_j^+,\ \  i\in\overline{1,M},\ (\alpha,j)\not\in I_0\cup I_1\cup ... \cup I_M,\\
|I_0|=d: &\ \   c_j^-<\left<e_j,Ay_{*,0}\right>< c_j^+,\ \  (\alpha,j)\not\in I_0.
\end{array}
\end{equation}
 
\noindent Specifically, if (\ref{blueinclusion}) and (\ref{y*ifeasible}) hold then each solution of (\ref{sp}) converges in finite time to (\ref{Fdef}) and, accordingly, the stress-vector $s(t)$ of $(D,A,C,R,l(t))$ converges in finite time to (\ref{SAF}).

\vskip0.2cm

\noindent {\bf Terminal response force.} The knowledge of the terminal value of stress-vector $s(t)$ allows to compute the {\it terminal response force} of elastoplastic system $(D,A,C,R,l(t))$. According to \cite[formulas (7) and (30)]{ESAIM}, the terminal response force of 
$(D,A,C,R,l(t))$ is the solution $r$ of
$$
  -D^Ts_*-D^TRr=0,\quad {\rm where}\quad s_*\in S.
$$

\subsection{Computations for the benchmark example}

\noindent Following \cite{SICON}, computation of (\ref{blueinclusion}) for the example of Fig.~\ref{fig1} gives 
$$
  \left(\hskip-0.1cm\begin{array}{c} 1\\ 0\\ 0\end{array}\hskip-0.1cm\right)\hskip-0.05cm\in{\rm cone}\left(\hskip-0.1cm\left(\begin{array}{ccccc}
  1 & 0 & 1 & 0 & 1\\
  0 & 0 & 1 & -1 & 1\\
  1 & -1 & 1 & 0 & 0
  \end{array}\right)\hskip-0.1cm\left\{\alpha e_j:(\alpha,j)\in I_0\right\}\hskip-0.1cm\right).
$$
One can observe that the only 4 possibilities for irreducible and admissible $I_0$ are
\begin{itemize}
    \item $I_0=\{(+,1),(-,3),(+,5)\}$, \item $I_0=\{(+,2),(+,3),(+,4)\}$,
    \item $I_0=\{(+,1),(+,2)\},$ 
    \item $I_0=\{(+,4),(+,5)\}$,
\end{itemize}
where, for clarity, we use "+" for $\alpha=1$ and "$-$" for $\alpha=-1$. 
The paper \cite{SICON} shows that feasibility condition (\ref{y*ifeasible})
for first two bullet points to take place are (\ref{FF2}) and (\ref{FF4}) respectively, and the corresponding terminal response force is given by formulas (\ref{FF1}) and (\ref{FF3}) accordingly. 

\end{document}